\documentclass[11pt,a4paper]{amsart}

\usepackage[T1]{fontenc}
\usepackage{a4wide}
\usepackage{mathtools,amssymb}
\usepackage[mathscr]{eucal}
\usepackage{amsthm}
\usepackage{enumitem}
\usepackage{microtype}
\usepackage{url}
\usepackage[colorlinks=true,linkcolor=blue,citecolor=blue,urlcolor=blue]{hyperref}

\numberwithin{equation}{section}
\newtheoremstyle{paperplain}%
  {6pt plus 2pt minus 2pt}
  {6pt plus 2pt minus 2pt}
  {\itshape}
  {}
  {\bfseries\upshape}
  {.}
  {0.5em}
  {}
\theoremstyle{paperplain}
\newtheorem{theorem}{Theorem}

\newtheorem{proposition}{Proposition}[section]
\theoremstyle{remark}

\newcommand{\Id}{\mathrm{id}}
\newcommand{\K}{\EuScript{K}}
\newcommand{\B}{\EuScript{B}}

\title[Hamana's injective envelope as a maximal rigid multiplier cover]{Hamana's injective envelope as\\ a maximal rigid multiplier cover}

\subjclass[2020]{Primary 46L05; Secondary 46L30, 54C55}
\keywords{Injective envelope, multiplier algebra, rigid extension, Gleason cover}
\date{}

\author[T.~Kania]{Tomasz Kania}
\address[T.~Kania]{Mathematical Institute\\Czech Academy of Sciences\\\v Zitn\'a 25 \\115 67 Praha 1\\Czech Republic  and  Institute of Mathematics and Computer Science\\ Jagiellonian University\\ {\L}ojasiewicza 6, 30-348 Krak\'{o}w, Poland
}
\email{kania@math.cas.cz, tomasz.marcin.kania@gmail.com}
\thanks{RVO: 67985840.}

\dedicatory{Dedicated to Aleksander B{\l}aszczyk on the occasion of his 79th birthday.}

\begin{document}

\begin{abstract}
Let $A$ be a unital $C^*$-algebra. An $A$-multiplier cover is a $C^*$-algebra $E$ together with a faithful non-degenerate $*$-homomorphism from $A$ to $M(E)$. We preorder such covers by $A$-preserving unital completely positive maps between their multiplier algebras. We prove that Hamana's injective envelope $I(A)$ is a greatest cover in this preorder and that the maximal rigid covers are precisely those whose multiplier algebra is canonically $*$-isomorphic to $I(A)$. Consequently, a maximal rigid cover is greatest, rather than merely maximal among rigid covers. For $A=C(X)$, we further classify these covers: after the canonical identification with $C(G(X))$, where $G(X)$ is the Gleason cover, their underlying ideals are exactly the algebras $C_0(U)$ for dense open $C^*$-embedded subsets $U$ of $G(X)$. Dense cozero subsets provide an important special case.
\end{abstract}

\maketitle

\section{Introduction}

Hamana's injective envelope $I(A)$ of a unital $C^*$-algebra $A$ is a canonical injective operator system containing $A$ and is itself a unital $C^*$-algebra \cite{Hamana1979OS,Hamana1979Cstar}. This note gives an order-theoretic characterisation of $I(A)$ in terms of multiplier algebras. The commutative motivation comes from B{\l}aszczyk's economical construction of the Gleason cover: one first maximises irreducibility and only then compactifies \cite{Blaszczyk}.

We collect all terminology here. A $*$-homomorphism $\iota\colon A\to M(E)$ is \emph{non-degenerate} when
\[
\overline{\operatorname{span}}\,\iota(A)E=E.
\]
Since $A$ is unital, non-degeneracy is equivalent to $\iota(1_A)=1_{M(E)}$: one implication follows because $\iota(1_A)$ acts as the identity on the dense subspace $\iota(A)E$, and the converse is immediate. An \emph{$A$-multiplier cover} is a pair $(E,\iota)$ in which $E$ is a $C^*$-algebra and $\iota\colon A\to M(E)$ is a faithful, equivalently injective, non-degenerate $*$-homomorphism. Thus $\iota$ is a unital $*$-monomorphism.

For two $A$-multiplier covers, write
\[
(E_1,\iota_1)\preccurlyeq(E_2,\iota_2)
\]
when there is a unital completely positive map $\Theta\colon M(E_1)\to M(E_2)$ such that $\Theta\circ\iota_1=\iota_2$. This is a preorder, but it need not be anti-symmetric. Indeed, for an infinite-dimensional Hilbert space $H$, the covers
\[
(\mathbb C,\Id_{\mathbb C})
\quad\text{and}\quad
(\K(H),\lambda\mapsto\lambda 1_{\B(H)})
\]
are comparable in both directions: scalar inclusion gives one comparison, and any state on $\B(H)$ gives the other.

A cover $(E,\iota)$ is \emph{rigid} if every unital completely positive map $\Phi\colon M(E)\to M(E)$ satisfying $\Phi\circ\iota=\iota$ equals $\Id_{M(E)}$. A rigid cover is \emph{maximal rigid} if
\[
(E,\iota)\preccurlyeq(F,\kappa),\qquad (F,\kappa)\ \text{rigid},
\]
implies $(F,\kappa)\preccurlyeq(E,\iota)$. On rigid covers, two-sided comparability is already canonical: the two witnessing maps have compositions fixing $A$, so rigidity makes them inverse complete order isomorphisms and hence inverse $*$-isomorphisms.

Let $j\colon A\to I(A)$ denote the canonical embedding. Our main result is the following.

\begin{theorem}\label{thm:A}
Let $A$ be a unital $C^*$-algebra and let $(E,\iota)$ be an $A$-multiplier cover. The following are equivalent.
\begin{enumerate}[label=\textup{(\roman*)},leftmargin=2.5em]
\item $(E,\iota)$ is maximal rigid;
\item $(E,\iota)$ is rigid and $(I(A),j)\preccurlyeq(E,\iota)$;
\item there is a unique $*$-isomorphism
\[
\Psi\colon M(E)\xrightarrow{\ \cong\ }I(A)
\]
such that $\Psi\circ\iota=j$.
\end{enumerate}
Whenever these conditions hold, $(E,\iota)$ is a greatest $A$-multiplier cover: $(F,\kappa)\preccurlyeq(E,\iota)$ for every $A$-multiplier cover $(F,\kappa)$.
\end{theorem}

Thus maximal rigid covers do not give a new construction of $I(A)$; Hamana's theorem is used in the proof. They do, however, characterise its multiplier-algebra realisations and show that maximality among rigid covers automatically upgrades to a greatest-element property.

For the commutative form, let $q\colon G(X)\twoheadrightarrow X$ be the Gleason cover of a compact Hausdorff space $X$. Recall that an open subset $U$ of a compact space $K$ is $C^*$-embedded when every bounded continuous complex-valued function on $U$ extends continuously to $K$.

\begin{theorem}\label{thm:B}
Let $X$ be compact Hausdorff.
\begin{enumerate}[label=\textup{(\roman*)},leftmargin=2.5em]
\item If $U\subseteq G(X)$ is dense, open, and $C^*$-embedded, then
\[
q_U^*\colon C(X)\longrightarrow M(C_0(U))=C_b(U),
\qquad q_U^*(f)=(f\circ q)|_U,
\]
makes $C_0(U)$ a maximal rigid $C(X)$-multiplier cover.
\item Conversely, if $(E,\iota)$ is a maximal rigid $C(X)$-multiplier cover, then the unique isomorphism
\[
\Psi\colon M(E)\xrightarrow{\ \cong\ }C(G(X))
\]
over $C(X)$ carries $E$ onto $C_0(U)$ for a dense open $C^*$-embedded subset $U$ of $G(X)$.
\end{enumerate}
In particular, every dense cozero subset of $G(X)$ gives a maximal rigid cover.
\end{theorem}

We work throughout with unital $C^*$-algebras. The non-unital case requires a different formulation.

\section{Maximal rigid multiplier covers}

We use the following properties of Hamana's injective envelope \cite{Hamana1979OS,Hamana1979Cstar}: $I(A)$ is an injective unital $C^*$-algebra, $j\colon A\to I(A)$ is a unital $*$-monomorphism, and every unital completely positive map $I(A)\to I(A)$ fixing $j(A)$ pointwise is the identity.

\begin{proposition}\label{prop:three-parts}
Let $A$ be a unital $C^*$-algebra.
\begin{enumerate}[label=\textup{(\roman*)},leftmargin=2.5em]
\item $(I(A),j)$ is a rigid $A$-multiplier cover.
\item For every $A$-multiplier cover $(E,\iota)$,
\[
(E,\iota)\preccurlyeq(I(A),j).
\]
Consequently, $(I(A),j)$ is a greatest $A$-multiplier cover and, in particular, is maximal rigid.
\item If $(E,\iota)$ is rigid and $(I(A),j)\preccurlyeq(E,\iota)$, there is a unique $*$-isomorphism $\Theta\colon M(E)\to I(A)$ satisfying $\Theta\circ\iota=j$.
\end{enumerate}
\end{proposition}

\begin{proof}
Since $I(A)$ is unital, $M(I(A))=I(A)$, and (i) is Hamana's rigidity.

For (ii), the unital $*$-isomorphism
\[
\iota(A)\longrightarrow j(A),\qquad \iota(a)\longmapsto j(a),
\]
extends, by injectivity of $I(A)$, to a unital completely positive map $\Theta\colon M(E)\to I(A)$ with $\Theta\circ\iota=j$. This proves the asserted comparison. If a rigid cover $(F,\kappa)$ satisfies $(I(A),j)\preccurlyeq(F,\kappa)$, part (iii) below gives the reverse comparison, so $(I(A),j)$ is maximal rigid.

For (iii), choose a unital completely positive map $\Lambda\colon I(A)\to M(E)$ with $\Lambda\circ j=\iota$, and let $\Theta\colon M(E)\to I(A)$ be supplied by (ii). The two compositions fix the respective copies of $A$, so rigidity gives
\[
\Theta\circ\Lambda=\Id_{I(A)}
\quad\text{and}\quad
\Lambda\circ\Theta=\Id_{M(E)}.
\]
Thus $\Theta$ and $\Lambda$ are inverse complete order isomorphisms. For completeness, this already forces $\Theta$ to be multiplicative. Indeed, for $x\in M(E)$, the Schwarz inequalities for $\Theta$ and $\Lambda$ give
\[
0\leqslant \Lambda\bigl(\Theta(x^*x)-\Theta(x)^*\Theta(x)\bigr)
=x^*x-\Lambda\bigl(\Theta(x)^*\Theta(x)\bigr)\leqslant 0.
\]
Since $\Lambda$ is injective, $\Theta(x^*x)=\Theta(x)^*\Theta(x)$; applying the same argument to $x^*$ shows that every $x$ lies in the multiplicative domain of $\Theta$. Hence $\Theta$ is a $*$-isomorphism.

If $\Theta'$ is another such $*$-isomorphism, then $\Lambda\circ\Theta'$ is a unital completely positive endomorphism of $M(E)$ fixing $\iota(A)$ pointwise. Rigidity yields $\Lambda\circ\Theta'=\Id_{M(E)}$, and therefore $\Theta'=\Theta$.
\end{proof}

\begin{proof}[Proof of Theorem~\ref{thm:A}]
Assume (i). Proposition~\ref{prop:three-parts}(ii) gives $(E,\iota)\preccurlyeq(I(A),j)$. Since $(I(A),j)$ is rigid, maximal rigidity of $(E,\iota)$ gives $(I(A),j)\preccurlyeq(E,\iota)$, proving (ii). The implication (ii)$\Rightarrow$(iii) is Proposition~\ref{prop:three-parts}(iii).

Assume (iii). Rigidity of $(I(A),j)$ transfers through $\Psi$, so $(E,\iota)$ is rigid. Moreover, Proposition~\ref{prop:three-parts}(ii) gives $(F,\kappa)\preccurlyeq(I(A),j)$ for every $A$-multiplier cover $(F,\kappa)$; composing a witnessing map with $\Psi^{-1}$ gives $(F,\kappa)\preccurlyeq(E,\iota)$. Thus $(E,\iota)$ is greatest and hence maximal rigid.
\end{proof}

\section{The commutative case}

B{\l}aszczyk \cite{Blaszczyk} constructed the Gleason cover by first passing to a maximal regular refinement for which the identity map onto the original compact space is irreducible, and then taking the \v{C}ech--Stone compactification. On the operator-algebraic side,
\[
I(C(X))\cong C(G(X))
\]
canonically over $C(X)$, where the embedding is $q^*f=f\circ q$ \cite{Hamana1979Cstar}.

\begin{proof}[Proof of Theorem~\ref{thm:B}]
Let $U\subseteq G(X)$ be dense, open, and $C^*$-embedded. Restriction gives an injective $*$-homomorphism $C(G(X))\to C_b(U)$, and $C^*$-embeddedness makes it surjective. Therefore
\[
M(C_0(U))=C_b(U)\cong C(G(X))\cong I(C(X))
\]
over $C(X)$. Theorem~\ref{thm:A} proves (i).

Conversely, let $(E,\iota)$ be maximal rigid. By Theorem~\ref{thm:A}, there is a unique $*$-isomorphism
\[
\Psi\colon M(E)\longrightarrow C(G(X))
\]
with $\Psi\circ\iota=q^*$. The image $J=\Psi(E)$ is a closed essential ideal of $C(G(X))$, so $J=C_0(U)$ for a unique dense open subset $U\subseteq G(X)$. The restriction $\Psi|_E\colon E\to C_0(U)$ extends uniquely to a $*$-isomorphism
\[
\widehat\Psi\colon M(E)\longrightarrow M(C_0(U))=C_b(U).
\]
If $\rho\colon C(G(X))\to C_b(U)$ denotes restriction, then $\widehat\Psi$ and $\rho\circ\Psi$ agree on the essential ideal $E$, and hence agree on $M(E)$. Since both $\Psi$ and $\widehat\Psi$ are surjective, $\rho$ is surjective. Thus every bounded continuous function on $U$ extends to $G(X)$, so $U$ is $C^*$-embedded. This proves (ii).

Finally, $G(X)$ is extremally disconnected and hence an $F$-space. Every cozero subset of an $F$-space is $C^*$-embedded \cite[Theorem~14.25]{GillmanJerison}, which gives the final assertion.
\end{proof}

\subsection*{Acknowledgements}
The author is grateful to the anonymous referee for a careful reading and for suggestions that improved both the organisation and the mathematical formulation of the paper. The author gratefully acknowledges support received from NCN Sonata-Bis~13 (2023/50/E/ST1/00067).

\paragraph{Conflict of Interest.}
The author declares that there is no conflict of interest.

\paragraph{Data Availability.}
No datasets were generated or analysed during the current study.

\end{document}